# Multi-state Operating Reserve Model of Aggregate Thermostatically-Controlled-Loads for Power System Short-Term Reliability Evaluation

Yi Ding[a], Wenqi Cui[a]*, Shujun Zhang[b], Hongxun Hui[a], Yiwei Qiu[c], Yonghua Song[a,d]

*a. College of Electrical Engineering, Zhejiang University, Hangzhou, China*
*b. State Grid Zhejiang Electric Power Company, Hangzhou, China*
*c. Department of Electrical Engineering, Tsinghua University, Beijing, China*
*d. Department of Electrical and Computer Engineering, University of Macau, Macau, China*

**Abstract**. Thermostatically-controlled-loads (TCLs) have been regarded as a good candidate for maintaining the power system reliability by providing operating reserve. The short-term reliability evaluation of power systems, which is essential for power system operators in decision making to secure the system real time balancing, calls for the accurate modelling of operating reserve provided by TCLs. However, the particular characteristics of TCLs make their dynamic response different from the traditional generating units, resulting in difficulties to accurately represent the reliability of operating reserve provided by TCLs with conventional reliability model. This paper proposes a novel multi-state reliability model of operating reserve provided by TCLs considering their dynamic response during the reserve deployment process. An analytical model for characterizing dynamics of operating reserve provided by TCLs is firstly developed based on the migration of TCLs' room temperature within the temperature hysteresis band. Then, considering the stochastic consumers' behaviour and ambient temperature, the probability distribution functions of reserve capacity provided by TCLs are obtained by cumulants. On this basis, the states of reserve capacity and the corresponding probabilities at each time instant are obtained for representing the reliability of operating reserve provided by TCLs in the $L_Z$-transform approach. Case studies are conducted to validate the proposed technique.

---

* Corresponding author: cuiwenqi@zju.edu.cn

Email Addresses: yiding@zju.edu.cn (Yi Ding), bluemoon_987@126.com (Shujun Zhang), huihongxun@zju.edu.cn (Hongxun Hui), ywqiu@mail.tsinghua.edu.cn (Yiwei Qiu), yhsong@umac.mo (Yonghua Song)

**Nomenclature**

*Acronyms*

| | |
|---|---|
| TCLs | Thermostatically-controlled-loads |
| N-TCLs | Thermostatically-controlled-loads in the ON mode |
| F-TCLs | Thermostatically-controlled-loads in the OFF mode |
| ORT | Operating reserve provided by TCLs |

*Variables and parameters*

| | |
|---|---|
| $\upsilon$ | Index of an individual TCL |
| $g$ | Index of a TCL group |
| $t$ | Index of time |
| $i$ | Index of a bus |
| $t_s$ | Deployment time of operating reserve provided by TCLs |
| $\xi$ | Index of the time interval when calculating dynamic aggregate response of TCLs |
| $\theta_\upsilon(t)$ | Room temperature corresponding to the $\upsilon$-th TCL at the time $t$ |
| $\theta_a(t)$ | Ambient temperature at the time $t$ |
| $\theta_{set,\upsilon}(t)$ | Set point temperature of the $\upsilon$-th TCL at the time $t$ |
| $\theta_{+,\upsilon}, \theta_{-,\upsilon}$ | Upper and lower temperature hysteresis band of the $\upsilon$-th TCL |
| $P_g(t)$ | Aggregate power of TCLs in group $g$ |
| $\Delta P_g(t)$ | Total deviation of TCLs' aggregate power in group $g$ |
| $T_{on,k_c}(t)$ | Expected on time of TCLs in the $c$-th cluster |
| $T_{off,k_c}(t)$ | Expected off time of TCLs in the $c$-th cluster |
| $\langle\tau\rangle\big|_{\theta_1}^{\theta_2}$ | Time duration $\tau$ for a TCL's room temperature to migrate from $\theta_1$ to $\theta_2$ |
| $f_{\Delta P_g(t)}(x)$ | Probability distribution function of $\Delta P_g(t)$ |
| $F_{P_g(t)}(x)$ | Cumulative distribution function of $P_g(t)$ |
| $RC_{i,j_i^{TCL}}$ | Reserve capacity of ORT at bus $i$ in the $j_i^{TCL}$ state |
| $\rho_{i,j_i^{TCL}}(t)$ | Probability of the reserve capacity $RC_{i,j_i^{TCL}}$ at bus $i$ for the time $t$ |
| $\text{MORT}_i(t)$ | Multi-state operating reserve provided by TCLs at bus $i$ for the time $t$ |
| $\text{MHOR}_i(t)$ | Multi-state hybrid operating reserve providers at bus $i$ for the time $t$ |
| $\text{MHGU}_i(t)$ | Multi-state hybrid generation units at bus $i$ for the time $t$ |
| $\text{MHGR}_i(t)$ | Multi-state hybrid generation and operating reserve provider at bus $i$ for the time $t$ |

# 1. Introduction

With the growing penetration of renewable energy sources (RES) in electric power system, both the power supply and demand have become highly time-varying, calling for a huge amount of balancing services to maintain the system reliability [1], [2]. Apart from increasing generation units, demand side resources (DSRs) with flexible and fast response capabilities have been regarded as an effective tool to enhance the system reliability by providing operating reserve [3], [4]. Because of the stochastic characteristic of RES and DSR, the electric power entities need to adequately aware of system operating pressures during a short interval, which could be achieved by power system short-term reliability evaluation [5].

Among all the DSRs, thermostatically-controlled loads (TCLs) have been considered as one of the most suitable resources for providing reserve services, since TCLs may account for over 40% of energy consumption for days with extremely hot or cold weather and can be flexibly controlled by adjusting thermostat set-point temperature [6]. It has been shown that 37% of load reduction in the confirmed demand response (DR) events in PJM in 2017 is provided by TCLs [7]. Such large share of participation makes it crucial to involve the dynamic response model of TCLs into the short-term reliability framework so that the system reliability level can be adequately assessed. However, one important characteristic of TCLs is that they cycle on and off by turns instead of consuming power at a constant value [8]. This may lead to the demand response rebound during the dynamic response process of TCLs [9], which is different from conventional generating units and therefore bring complexities in reliability evaluation and management [10].

References that have addressed the impact of DR on power system reliabilities, such as [11], [12], usually consider the general steady state model but not the accurate dynamic response model. Most previous researches represent the demand response by two-state of demand reduction and its failure, which is the same form of conventional reliability model for generation unit [13]. Reference [12] extracts the demand reduction and failure possibilities of different types of DSRs from the historical data gathered through surveys. The potential impacts of DR on reliability of a residential distribution network is quantified and verified, respectively. Reference [14] extents the two-state model to a multi-state model by involving several derated states during the DR process. State transition diagram is applied to represent the transition mode from one state to another. However, such general model for DSRs is difficult to represent the dynamic aggregate response of TCLs, which is different from conventional generation units because of the cyclical operation characteristics of TCLs. Moreover, considering that aggregate response of DSRs is influenced by stochastic consumers' behaviors, environment parameters and control algorithms, it may be unpractical to obtain the time-varying state transition diagram at each time instant.

In order to capture the dynamics of power system components (e.g., wind farms, coal thermal generators, etc.), multi-state models are increasingly adopted in power system reliability evaluation [14]. Instead of the conventional two-state model that only involves steady state and its failure, multi-state models achieve higher accuracy by involving in several states and probabilities according to the actual output distribution of the components [15], [16]. The complex dynamic response characteristic of aggregate TCLs means that it is also much more rational to model ORT as multi-state operating reserve provider instead of conventional two-state ones [14]. In order to represent the increased states and probabilities brought by multi-state model, the Universal Generating Functions (UGF) is proposed in [17], [18] to algebraically find the entire multi-state system performance distributions through the steady-states performance distributions of its elements. The effectiveness of UGF in power system long-term reliability evaluation has been verified in [19]. As an extension of conventional UGF, $L_z$-transform approach is put forward to involve in the time-varying probabilities of different states so that the dynamic reliability of the multi-state system can be evaluated [20], [21]. This allows the $L_z$-transform approach to be applied to power system short-term reliability evaluation considering hybrid generation and reserve providers [10]. The precondition for the application of $L_z$-transform approach in the reliability modelling of the operating reserve is to obtain the states of reserve capacity and the corresponding probabilities. However, there is a lack of the multi-state model for the operating reserve provided by TCLs, making it difficult to integrate with the multi-state power system reliability evaluation techniques to achieve an accurate assessment of the system operating pressures during a short interval.

The major gap in achieving the multi-state operating reserve model of TCLs is the lack of methods to directly obtain the dynamic aggregate power of TCLs controlled for reserve deployment under uncertainties. Existing methods for modeling the dynamic response of TCLs can be classified into two categories, Monte Carlo simulation approaches [22] and analytical techniques [23], [24]. Monte Carlo simulation can accurately model the influence of multiple stochastic factors but is time-consuming, which is difficult to fulfill the requirements of short-term reliability evaluation on computational time [22]. Existing analytical models, including the state-based model, the job scheduling model, and Gray-box model, focus on the relationship between the control signal and the temperature density evolution of TCLs based on iteration [24]. In state-based models, variations of TCLs' room temperature are classified into several state bins and represented by a state matrix [25], which is widely adopted in designing optimal control strategies for TCLs [26], [27]. However, TCLs are assumed to migrate with the same speed in the same temperature state bins [23], resulting in the difficulties in applying to large populations of TCLs with widely distributed parameters [28], [24]. With the aim of finding the optimal demand reduction, regulation of TCLs in [29] is converted into the job scheduling problem solved by greedy algorithm and binary search algorithm. Gray-box model is proposed in [30] to represent TCLs integrated with ON/OFF controllers using the data-driven technique. Although these models can capture the dynamic response of TCLs, the aggregate power cannot be directly obtained from these models. This results in the difficulties to obtain the states of reserve capacity and the corresponding time-varying probabilities for representing operating reserve provided by TCLs in a multi-state manner.

This paper proposes a novel multi-state reliability model of operating reserve provided by TCLs considering their dynamic response characteristics. The provision of operating reserve by TCLs is modelled as a discrete-state continuous-time process, which is represented by $L_Z$-transform approach proposed in [20]. Firstly, the dynamic aggregate power of TCLs is generated directly from the migration of TCLs' room temperature during the reserve deployment process. Then, considering the stochastic consumers' behavior and ambient temperature, the property of cumulants is applied to obtain the probability distribution functions of reserve capacity provided by TCLs. On this basis, the states of reserve capacity and the corresponding probabilities at each time instant are obtained for representing the reliability of ORT in the $L_Z$-transform approach. In this way, the system reliability with hybrid generation and operating reserve providers are also developed. The IEEE Reliability Test System (RTS) is applied to illustrate the validity and benefits of the proposed technique [10]. The major contributions of this paper are as follows:

1) A novel analytical model to characterize the dynamic response of aggregate TCLs for the provision of operating reserve, which can directly obtain the dynamic aggregate power of TCLs without the time-consuming iteration, is proposed.
2) A probabilistic model of operating reserve provided by TCLs, which considers the stochastic consumers' behavior and the ambient temperature, is presented to obtain the time varying probability distribution of reserve capacity.
3) A multi-state reliability model of operating reserve provided by TCLs, which reflects TCLs' particular dynamics and characteristics (e.g., the demand response rebound), is developed for the power system short-term reliability evaluation based on $L_Z$-transform approach.

## 2. Analysis of TCLs' Dynamic Response for the Provision of Operating Reserve

*2.1. Framework to Obtain the Multi-State Model of Operating Reserve Provided by TCLs*

Operating reserve (OR), which consists of spinning and non-spinning reserve according to the definition by the North American Electric Reliability Corporation (NERC), refers to the stand-by power or demand reduction that can be called on with short notice to deal with an unexpected mismatch between generation and load [31], [32]. The power consumption of TCLs can be easily controlled by changing set point temperature with a short time, making it suitable to provide operating reserve by reducing power consumption according to the reserve deployment instructions [33], [34]. Considering that the power consumption of an individual TCL is too small, TCLs are usually aggregated as an equivalent operating reserve provider (ORT) to provide operating reserve to the power system [35], [36]. Then, the system operator can dispatch the TCLs as the traditional operating reserve to enhance the reliability of the power systems [37].

This paper proposes a novel reliability model of ORT for the power system short-term reliability evaluation considering the dynamic response of TCLs, as illustrated in Fig. 1. Firstly, an analytical model to characterize the dynamics of TCLs controlled for providing ORT is developed based on the model of an individual TCL. Secondly, a stochastic model for ORT is proposed based on the dynamic response of TCLs considering the stochastic consumers' behavior and the ambient temperature. Thirdly, $L_Z$-transform approach is applied to construct the multi-state reliability model of ORT. These models will be elaborated in the following parts of the paper.

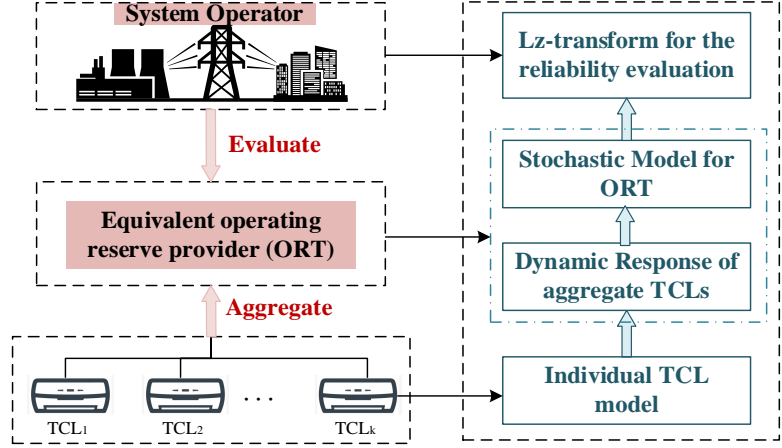

Fig. 1. Framework to obtain the multi-state model of operating reserve provided by TCLs

## 2.2. Equivalent Model of Operating Reserve Provided by TCLs

The operation process of an individual TCL is described as the following hybrid state model [38]:

$$\frac{d\theta_\upsilon(t)}{dt} = -\frac{1}{C_\upsilon R_\upsilon}[\theta_\upsilon(t) - \theta_a(t) + m_\upsilon(t)R_\upsilon Q_\upsilon] \tag{1}$$

$$m_\upsilon(t) = \begin{cases} 1, & \theta_\upsilon(t) > \theta_{+,\upsilon} \\ 0, & \theta_\upsilon(t) < \theta_{-,\upsilon} \\ m_\upsilon(t-1), & otherwise \end{cases} \tag{2}$$

where $\theta_\upsilon(t)$ is the room temperature corresponding to the $\upsilon$-th TCL at time $t$, $\theta_a(t)$ is the ambient temperature. $C_\upsilon$ and $R_\upsilon$ are the thermal capacity and thermal resistance corresponding to the room of the $\upsilon$-th TCL, respectively. $Q_\upsilon$ is the energy transfer rate of the $\upsilon$-th TCL, which is equal to the product of the input power $p_\upsilon$ and the coefficient of performance $COP_\upsilon$ of the $\upsilon$-th TCL. $m_\upsilon(t)$ represents the ON or OFF mode of the $\upsilon$-th TCL. $\theta_\upsilon(t)$ is maintained around its set point temperature $\theta_{set,\upsilon}$ with a dead band of $\Delta\theta_\upsilon$ by switching *on* ($m_\upsilon(t)=1$) or switching *off* ($m_\upsilon(t)=0$) TCL compressor. The temperature range between the lower band ($\theta_{-,\upsilon}^0 = \theta_{set,\upsilon} - 0.5 \times \Delta\theta_\upsilon$) and the upper band ($\theta_{+,\upsilon}^0 = \theta_{set,\upsilon} + 0.5 \times \Delta\theta_\upsilon$) is defined as the temperature hysteresis band [$\theta_{-,\upsilon}^0$, $\theta_{+,\upsilon}^0$]. Eq. (2) illustrates the changes of $m_\upsilon(t)$ corresponding to the temperature hysteresis band when the TCL operates for cooling in summer.

A typical curve of the consumed power of an individual TCL and the corresponding variation of room temperature is illustrated in Fig. 2(a), where the set point temperature is increased by $\beta_\upsilon$ at the time $t_s$. It can be observed from Fig. 2(a) that the TCL operates cyclically within the temperature hysteresis band. Changes of the set point temperature will result in the changes of the time that TCLs spend in the ON mode and the OFF mode, thereby changes the consumed power a TCL [39].

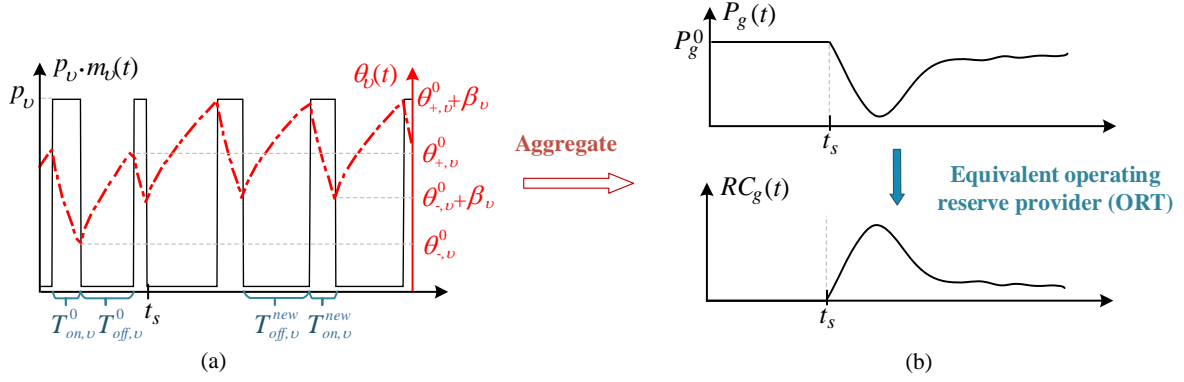

Fig. 2. Equivalent operating reserve provided by aggregate TCLs (a) The consumed power and the corresponding variation of room temperature of an individual TCL (b) Aggregate power and equivalent operating reserve of TCLs in group $g$ after the changes of set point temperature

TCLs provide operating reserve by actively reducing aggregate power consumption through changing set point temperature [40]. Fig. 2(b) is the aggregate power $P_g(t)$ of a TCL group $g$ after the changes of set point temperature at the time $t_s$. $P_g^0$ denotes the initial aggregate power of TCLs before $t_s$. Then, as illustrated by Fig. 2(b), reserve capacity $RC_g(t)$ of operating reserve provided by TCLs in the group $g$ is the difference of aggregate power before and after the changes of set point temperature [41]:

$$RC_g(t) = P_g^0 - P_g(t) \tag{3}$$

Therefore, in order to know the exact value of reserve capacity, the key point is to obtain the value of $P_g(t)$ after the changes of set point temperature. The analytical model of $P_g(t)$ will be discussed in detail in Section 3. Then, the variation of $P_g(t)$ considering the stochastic changes of set point temperature and the ambient temperature is generated in Section 4. On this basis, the multi-state model of operating reserve provided by TCLs can be obtained.

## 3. Analytical Model to Characterize the Dynamic Response of Heterogeneous TCLs for the Provision of Operating Reserve

### 3.1. General Representation of TCLs' Aggregate Power

Eqs. (1)-(2) illustrate that TCLs operate cyclically within the temperature hysteresis band. The expected duration of the $\upsilon$-th TCL in the ON mode and OFF mode at the time $t$ are denoted by $T_{on,\upsilon}(t)$ and $T_{off,\upsilon}(t)$, respectively. Considering the cyclical operation characteristic of TCLs, the average power $\bar{p}_\upsilon$ of the $\upsilon$-th TCL can be obtained by the ratio of $T_{on,\upsilon}(t)$ to the whole duty cycle ($T_{on,\upsilon}(t)+T_{off,\upsilon}(t)$).

Denoting $\Gamma$ as the set of all the TCLs in group $g$, the aggregate power $P_g(t)$ of group $g$ is the summation of $\bar{p}_\upsilon$:

$$P_g(t) = \sum_{\upsilon \in \Gamma} \bar{p}_\upsilon(t) = \sum_{\upsilon \in \Gamma} p_\upsilon \cdot \frac{T_{on,\upsilon}(t)}{T_{on,\upsilon}(t) + T_{off,\upsilon}(t)} \tag{4}$$

For a group of homogeneous TCLs, the parameters of TCLs are similar, leading to similar on time $T_{on,\upsilon}(t)$ and off time $T_{off,\upsilon}(t)$ of each individual TCL in this group. $g_c$ denotes the typical TCL in the group $g$ with homogeneous TCLs, then (4) can be approximated by:

$$P_g(t) = \frac{T_{on,g_c}(t)}{T_{on,g_c}(t) + T_{off,g_c}(t)} \cdot \sum_{\upsilon \in \Gamma} p_\upsilon \tag{5}$$

In this way, the aggregate power of group $g$ is converted to the calculation of $T_{on,g_c}(t)$ and $T_{off,g_c}(t)$. To deal with load heterogeneity, all the TCLs is classified into $Q$ clusters according to the on time and off time of each TCL using the k-means algorithm [42]. Let $k_c$ denotes the center of the $c$-th cluster and corresponds to the $k_c$-th TCL in group $g$; $S_c$ denotes the set of all the ACs belong to the $c$-th cluster. The total aggregate power of TCLs in $\Gamma$ is approximated by the summation of TCLs' aggregate power in each cluster:

$$P_g(t) \approx \sum_{c=1}^{Q} \left( \frac{T_{on,k_c}(t)}{T_{on,k_c}(t) + T_{off,k_c}(t)} \cdot \sum_{\upsilon \in S_c} p_\upsilon \right) \tag{6}$$

Changes of hysteresis band will result in the changes of $T_{on,k_c}(t)$ and $T_{off,k_c}(t)$, thereby influence the level of aggregate power and thus provide operating reserve to power systems. Therefore, the key point to obtain the dynamic aggregate power of TCLs after the changes of set point temperature is to obtain the $T_{on,k_c}(t)$ and $T_{off,k_c}(t)$ at each time instant.

*3.2. Aggregate Response of Heterogeneous TCLs*

The expected on time $T_{on,k_c}(t)$ and off time $T_{off,k_c}(t)$ of the TCLs in the $c$-th cluster at each time instant is corresponding to the migration of TCLs within the hysteresis band. Existing research studies have found that sudden changes of the set point temperature will cause temporary synchronization of TCLs, resulting in large power fluctuations. Therefore, the safe protocol-2 introduced in [43] is utilized in this paper to avoid temporary synchronization of the TCLs. In this way, the migration of the TCLs in the $c$-th cluster after the shifting of temperature hysteresis band is shown in Fig. 3, where $\theta_L^{on/off}$ and $\theta_H^{on/off}$ are the lowest boundary and highest boundary corresponding to the room temperature of TCLs in the ON/OFF mode, respectively. $T_{on,k_c}(t)$ and $T_{off,k_c}(t)$ can be calculated according to TCLs' room temperature range. To help illustrate the migration of TCLs' room temperature range, the TCLs in the ON mode and OFF mode are abbreviated as N-TCLs and F-TCLs, respectively.

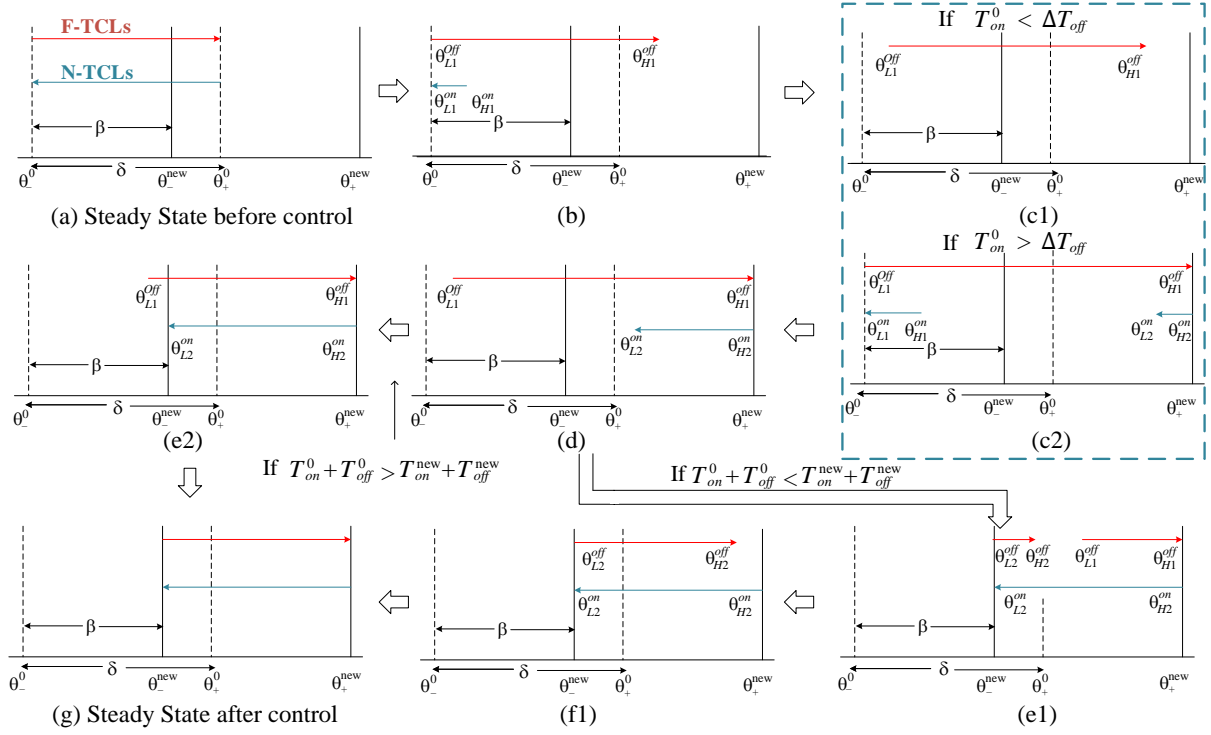

Fig. 3. Migration of TCLs' room temperature after the shifting of the temperature hysteresis band

Before the reserve deployment time $t_s$, TCLs operate cyclically within the initial temperature hysteresis band (Fig. 3(a)). The steady state on time $T^0_{on,k_c}$ and off time $T^0_{off,k_c}$ of the $c$-th cluster corresponding to the initial temperature hysteresis band [$\theta^0_{-,k_c}$, $\theta^0_{+,k_c}$] can be calculated from Eqs. (1)-(2) as follows [23]:

$$T^0_{on,k_c} = C_{k_c} R_{k_c} \ln(\frac{p_{k_c} R_{k_c} + \theta^0_{+,k_c} - \theta_a}{p_{k_c} R_{k_c} + \theta^0_{-,k_c} - \theta_a}) \tag{7}$$

$$T^0_{off,k_c} = C_{k_c} R_{k_c} \ln(\frac{\theta_a - \theta^0_{-,k_c}}{\theta_a - \theta^0_{+,k_c}}) \tag{8}$$

The temperature hysteresis band of the $k_c$-th TCL after the shifting of temperature hysteresis band by $\beta_{k_c}$ at the time $t_s$ is denoted as [$\theta^{new}_{-,k_c}, \theta^{new}_{+,k_c}$], which equals to [$\theta^0_{-,k_c}+\beta_{k_c}, \theta^0_{+,k_c}+\beta_{k_c}$]. The steady state cooling time $T^{new}_{on,k_c}$ and heating time $T^{new}_{off,k_c}$ of the $c$-th cluster corresponding to the new temperature hysteresis band are also obtained from Eqs. (1)-(2) as follows:

$$T^{new}_{on,k_c} = C_{k_c} R_{k_c} \ln(\frac{p_{k_c} R_{k_c} + \theta^0_{+,k_c} - \theta_a + \beta_{k_c}}{p_{k_c} R_{k_c} + \theta^0_{-,k_c} - \theta_a + \beta_{k_c}}) \tag{9}$$

$$T^{new}_{off,k_c} = C_{k_c} R_{k_c} \ln(\frac{\theta_a - \theta^0_{-,k_c} - \beta_{k_c}}{\theta_a - \theta^0_{+,k_c} - \beta_{k_c}}) \tag{10}$$

The duration $\Delta T_{off,k_c}$ for the $c$-th cluster to migrate from $\theta^0_{+,k_c}$ to $\theta^{new}_{+,k_c}$ is:

$$\Delta T_{off,k_c} = C_{k_c} R_{k_c} \ln(\frac{\theta_a - \theta^0_{+,k_c}}{\theta_a - \theta^0_{+,k_c} - \beta_{k_c}}) \tag{11}$$

When TCLs are migrating to the new temperature hysteresis band, the range of TCLs' room temperature in the ON mode and OFF mode will follow the process shown by Fig. 3(b) to Fig. 3(f). If it takes time duration $\tau$ for a TCL's room temperature to migrate from $\theta_1$ to $\theta_2$, $\tau$ is labeled as $\langle \tau \rangle \big|_{\theta_1}^{\theta_2}$ to indicate the time duration corresponding to $[\theta_1, \theta_2]$.

*3.3. Dynamics of TCLs During their Migration to the New Temperature Hysteresis Band*

As mentioned above in Eq. (6), the key point to obtain the dynamic aggregate power of TCLs after the changes of set point temperature is to obtain the expected on time $T_{on,k_c}(t)$ and expected off time $T_{off,k_c}(t)$ at each time instant. The value of $T_{on,k_c}(t)$ and $T_{off,k_c}(t)$ obtained from the migration of TCLs within the hysteresis band corresponding to the process shown by Fig. 3(b) to Fig. 3(f) are illustrated as follows.

*1) Migration of TCLs' room temperature shown by Fig. 3(b)*

It takes the time $\Delta T_{off,k_c}$ for the F-TCLs in Fig. 3(a) to reach $\theta_{+,k_c}^{new}$. Meanwhile, it takes $T_{on,k_c}^0$ for all the N-TCLs in Fig. 3(a) to switch from ON mode to OFF mode. Hence, if $t - t_s < \min\{\Delta T_{off,k_c}, T_{on,k_c}^0\}$, TCLs will migrate from Fig. 3(a) to Fig. 3(b). It takes F-TCLs in Fig. 3(a) the time ($t$-$t_s$) to migrate from $\theta_{+,k_c}^0$ to $\theta_{H1,k_c}^{off}$. Meanwhile, it takes N-TCLs in Fig. 3(a) the time ($t$-$t_s$) to migrate from $\theta_{+,k_c}^0$ to $\theta_{H1,k_c}^{off}$. Therefore, the expected on time $T_{on,k_c}(t)$ and expected off time $T_{off,k_c}(t)$ are:

$$\begin{cases} \text{If } t_s \leq t < t_s + \min\{\Delta T_{off,k_c}, T_{on,k_c}^0\} \\ T_{on,k_c}(t) = \left\langle T_{on,k_c}^0 - (t-t_s) \right\rangle \Big|_{\theta_{H1,k_c}^{On}}^{\theta_{L1,k_c}^{On}} = T_{on,k_c}^0 - (t-t_s) \\ T_{off,k_c}(t) = \left\langle T_{off,k_c}^0 + (t-t_s) \right\rangle \Big|_{\theta_{L1,k_c}^{Off}}^{\theta_{H1,k_c}^{Off}} = T_{off,k_c}^0 + (t-t_s) \end{cases} \quad (12)$$

*2) Migration of TCLs' room temperature shown by Fig. 3(c)*

- If $T_{on,k_c}^0 < t - t_s < \Delta T_{off,k_c}$, all the N-TCLs in Fig. 3(b) have switched from ON mode to OFF mode before the F-TCLs in Fig. 3(b) have reached $\theta_{+,k_c}^{new}$. Hence, TCLs will migrate from Fig. 3(b) to Fig. 3(c1). It takes the time $T_{on,k_c}^0$ for the last N-TCL to migrate from $\theta_{+,k_c}^0$ to $\theta_{-,k_c}^0$, after which the TCL become the last F-TCL and spend the time $(t-t_s) - T_{on,k_c}^0$ to migrate from $\theta_{-,k_c}^0$ to $\theta_{L1,k_c}^{off}$. Therefore, the expected on time $T_{on,k_c}(t)$ and expected off time $T_{off,k_c}(t)$ are:

$$\begin{cases} \text{If } T_{on,k_c}^0 + t_s < t < \Delta T_{off,k_c} + t_s \\ T_{on,k_c}(t) = 0 \\ T_{off,k_c}(t) = \left\langle T_{off,k_c}^0 + (t-t_s) \right\rangle \Big|_{\theta_{-,k_c}^0}^{\theta_{H1,k_c}^{Off}} - \left\langle (t-t_s) - T_{on,k_c}^0 \right\rangle \Big|_{\theta_{-,k_c}^0}^{\theta_{L1,k_c}^{Off}} = T_{off,k_c}^0 + T_{on,k_c}^0 \end{cases} \quad (13)$$

- If $\Delta T_{off,k_c} < t - t_s < T_{on,k_c}^0$, TCLs will migrate from Fig. 3(b) to Fig. 3(c2). It takes the time $\Delta T_{off,k_c}$ for the first F-TCL in Fig. 3 (b) to migrate from $\theta_{+,k_c}^0$ to $\theta_{+,k_c}^{new}$, after which the TCL will switch from OFF mode to ON mode and spend the time $(t-t_s)-\Delta T_{off}$ to migrate from $\theta_{+,k_c}^{new}$ to $\theta_{L2,k_c}^{on}$. Therefore, the expected on time $T_{on,k_c}(t)$ and expected off time $T_{off,k_c}(t)$ are:

$$\begin{cases} \text{If} \quad \Delta T_{off,k_c}+t_s < t < T_{on,k_c}^0+t_s \\ T_{on,k_c}=\left\langle T_{on,k_c}^0-(t-t_s)\right\rangle\Big|_{\theta_{H1,k_c}^{on}}^{\theta_{L1,k_c}^{on}}+\left\langle (t-t_s)-\Delta T_{off,k_c}\right\rangle\Big|_{\theta_{+,k_c}^{new}}^{\theta_{L2,k_c}^{on}}=T_{on,k_c}^0-\Delta T_{off,k_c} \\ T_{off,k_c}=\left\langle T_{off,k_c}^0+\Delta T_{off,k_c}\right\rangle\Big|_{\theta_{-,k_c}^0}^{\theta_{+,k_c}^{new}}=T_{off,k_c}^0+\Delta T_{off,k_c} \end{cases} \quad (14)$$

*3) Migration of TCLs' room temperature shown by Fig. 3(d)*

It takes the time ($\Delta T_{off,k_c}+T_{on,k_c}^{new}$) for TCLs in Fig. 3(c) to reach $\theta_{-,k_c}^{new}$. Therefore, if $\max\{\Delta T_{off,k_c}, T_{on,k_c}^0\}<t-t_s<\Delta T_{off,k_c}+T_{on,k_c}^{new}$, TCLs will migrate from Fig. 3(c) to Fig. 3(d). The expected on time $T_{on,k_c}(t)$ and expected off time $T_{off,k_c}(t)$ are:

$$\begin{cases} \text{If} \quad \max\{\Delta T_{off,k_c}, T_{on,k_c}^0\}+t_s < t < \Delta T_{off,k_c}+T_{on,k_c}^{new}+t_s \\ T_{on,k_c}(t)=\left\langle t-t_s-\Delta T_{off,k_c}\right\rangle\Big|_{\theta_{H2,k_c}^{on}}^{\theta_{L2,k_c}^{on}}=t-t_s-\Delta T_{of,k_cf} \\ T_{off,k_c}(t)=\left\langle T_{off,k_c}^0+\Delta T_{off,k_c}\right\rangle\Big|_{\theta_{-,k_c}^0}^{\theta_{+,k_c}^{new}}-\left\langle (t-t_s)-T_{on,k_c}^0\right\rangle\Big|_{\theta_{-,k_c}^0}^{\theta_{L1,k_c}^{off}}=T_{off,k_c}^0+T_{on,k_c}^0+\Delta T_{off,k_c}-t+t_s \end{cases} \quad (15)$$

*4) Migration of TCLs' room temperature shown by Fig. 3(e)*

It takes the time $(T_{on,k_c}^0+T_{off,k_c}^0-T_{off,k_c}^{new}-\Delta T_{off,k_c})$ for all the F-TCLs in Fig. 3(d) to reach $\theta_{-,k_c}^{new}$. Meanwhile, it takes the time ($\Delta T_{off,k_c}+T_{on,k_c}^{new}$) for the first N-TCL in Fig. 3(d) to reach $\theta_{-,k_c}^{new}$.

- If $T_{on,k_c}^0+T_{off,k_c}^0-(T_{off,k_c}^{new}-\Delta T_{off,k_c})<\Delta T_{off,k_c}+T_{on,k_c}^{new}$, namely, $T_{on,k_c}^0+T_{off,k_c}^0<T_{on,k_c}^{new}+T_{off,k_c}^{new}$, the last F-TCL in Fig. 3(d) reaches $\theta_{-,k_c}^{new}$ before the first N-TCL in Fig. 3(d). Hence, there would exist a gap during the migration of TCLs, as is shown in Fig. 3 (e1) and Fig. 3(f1). In this case, if $t-t_s>\Delta T_{off,k_c}+T_{on,k_c}^{new}$, after the first N-TCL in Fig. 3(d) reaches $\theta_{-,k_c}^{new}$, it will switch to OFF mode and spend the remaining time $((t-t_s)-(\Delta T_{off,k_c}+T_{on,k_c}^{new}))$ to migrate from $\theta_{L2,k_c}^{off}$ to $\theta_{H2,k_c}^{off}$. Meanwhile, if $t-t_s<T_{on,k_c}^0+T_{off,k_c}^0+\Delta T_{off,k_c}$, the last F-TCL in Fig. 3(d) have not reached $\theta_{+,k_c}^{new}$. Hence, TCLs will migrate from Fig. 3(d) to Fig. 3(e1). The expected on time $T_{on,k_c}(t)$ and expected off time $T_{off,k_c}(t)$ are:

$$\begin{cases} \text{If } \Delta T_{off,k_c}+T_{on,k_c}^{new}+t_s<t<T_{on,k_c}^0+T_{off,k_c}^0+\Delta T_{off,k_c}+t_s \ \& \ T_{on,k_c}^0+T_{off,k_c}^0<T_{on,k_c}^{new}+T_{off,k_c}^{new} \\ T_{on,k_c}(t)=\left\langle T_{on,k_c}^{new}\right\rangle\Big|_{\theta_{+,k_c}^{new}}^{\theta_{-,k_c}^{new}}=T_{on,k_c}^{new} \\ T_{off,k_c}(t)=\left\langle T_{off,k_c}^0+T_{on,k_c}^0+\Delta T_{off,k_c}-t+t_s\right\rangle\Big|_{\theta_{L1,k_c}^{off}}^{\theta_{H1,k_c}^{off}}+\left\langle (t-t_s)-(\Delta T_{off,k_c}+T_{on,k_c}^{new})\right\rangle\Big|_{\theta_{L2,k_c}^{off}}^{\theta_{H2,k_c}^{off}}=T_{off,k_c}^0+T_{on,k_c}^0-T_{on,k_c}^{new} \end{cases} \quad (16)$$

Then, if $t-t_s>T_{on,k_c}^0+T_{off,k_c}^0+\Delta T_{off,k_c}$, all the F-TCLs in [$\theta_{L1,k_c}^{off}$, $\theta_{H1,k_c}^{off}$] of Fig. 3(e1) have reached $\theta_{+,k_c}^{new}$ and switched to the ON mode. Meanwhile, if $t-t_s<T_{on,k_c}^{new}+T_{off,k_c}^{new}+\Delta T_{off,k_c}$, the F-TCLs corresponding to the temperature range [$\theta_{L2,k_c}^{off}$, $\theta_{H2,k_c}^{off}$] in Fig. 3(e1) have not reached $\theta_{+,k_c}^{new}$. Therefore, TCLs will migrate from Fig. 3(e1) to Fig. 3(f1). The expected on time $T_{on,k_c}(t)$ and expected off time $T_{off,k_c}(t)$ are:

$$\begin{cases} \text{If} \quad T_{on,k_c}^0 + T_{off,k_c}^0 + \Delta T_{off,k_c} + t_s < t < T_{on,k_c}^{new} + T_{off,k_c}^{new} + \Delta T_{off,k_c} + t_s \\ T_{on,k_c}(t) = \left\langle T_{on,k_c}^{new} \right\rangle \Big|_{\theta_+^{new}}^{\theta_-^{new}} = T_{on,k_c}^{new} \\ T_{off,k_c}(t) = \left\langle (t-t_s) - (\Delta T_{off,k_c} + T_{on,k_c}^{new}) \right\rangle \Big|_{\theta_{L2,k_c}^{off}}^{\theta_{H2,k_c}^{off}} = t - t_s - \Delta T_{off,k_c} - T_{on,k_c}^{new} \end{cases} \quad (17)$$

- If $\Delta T_{off,k_c} + T_{on,k_c}^{new} < t - t_s < T_{on,k_c}^0 + T_{off,k_c}^0 - (T_{off,k_c}^{new} - \Delta T_{off,k_c})$, TCLs will migrate from Fig. 3(d) to Fig. 3(e2). The expected on time $T_{on,k_c}(t)$ and expected off time $T_{off,k_c}(t)$ are:

$$\begin{cases} \text{If} \quad \Delta T_{off,k_c} + T_{on,k_c}^{new} + t_s < t < T_{on,k_c}^0 + T_{off,k_c}^0 - T_{off,k_c}^{new} + \Delta T_{off,k_c} + t_s \\ T_{on,k_c}(t) = T_{on,k_c}^{new} \\ T_{off,k_c}(t) = \left\langle T_{off,k_c}^0 + T_{on,k_c}^0 + \Delta T_{off,k_c} - t + t_s \right\rangle \Big|_{\theta_{L1,k_c}^{off}}^{\theta_{+,k_c}^{new}} = T_{off,k_c}^0 + T_{on,k_c}^0 + \Delta T_{off,k_c} - t + t_s \end{cases} \quad (18)$$

5) *Migration of TCLs' room temperature shown by Fig. 3(g)*

After the process illustrated by Fig. 3(f1) and Fig. 3(e2), all the TCLs will be covered in the new temperature hysteresis band. The expected on time $T_{on,k_c}(t)$ and expected off time $T_{off,k_c}(t)$ are:

$$\begin{cases} \text{If} \ t > T_{on,k_c}^{new} + T_{off,k_c}^{new} + \Delta T_{off,k_c} + t_s \ \& \ T_{on,k_c}^0 + T_{off,k_c}^0 < T_{on,k_c}^{new} + T_{off,k_c}^{new} \\ or \ \text{If} \ t > T_{on,k_c}^0 + T_{off,k_c}^0 - T_{off,k_c}^{new} + \Delta T_{off,k_c} + t_s \ \& \ T_{on,k_c}^0 + T_{off,k_c}^0 > T_{on,k_c}^{new} + T_{off,k_c}^{new} \\ T_{on,k_c}(t) = \left\langle T_{on,k_c}^{new} \right\rangle \Big|_{\theta_+^{new}}^{\theta_-^{new}} = T_{on,k_c}^{new} \\ T_{off,k_c}(t) = \left\langle T_{off,k_c}^{new} \right\rangle \Big|_{\theta_-^{new}}^{\theta_+^{new}} = T_{off,k_c}^{new} \end{cases} \quad (19)$$

## 4. Multi-state Reliability Model of Operating Reserve Provided by TCLs

This section obtains the probability distribution of TCLs' aggregate power affected by the variations of ambient temperature and set point temperature. On this basis, the multi-state model of operating reserve provided by TCLs is obtained by $L_z$-transform approach.

### 4.1. Probabilistic Model of Operating Reserve Provided by TCLs

In practice, it is observed that TCLs' power consumptions are correlated with the ambient temperature and the setpoint temperature [44]. Instead of the time-consuming Monte Carlo Simulation method, the property of cumulants is applied to compute the probability distribution of aggregate power in a systematic way [45].

It can be seen from Eqs. (12)-(19) that the calculation of aggregate power is corresponding to different time intervals. Because of the variation of ambient temperature $\Delta\theta_a(t)$ and the variation of set point temperature $\Delta\theta_{set,k_c}(t)$, the endpoints of each time interval are also uncertain. Hence, $T_{on,k_c}(t)$ and $T_{off,k_c}(t)$ are not deterministically determined by just one equation in Eqs.(12)-(19). $\Xi$ denotes the number of time interval corresponding to migration of TCLs shown by Fig. 3. The $T_{on,k_c}(t)$ and $T_{off,k_c}(t)$ in Eqs. (12)-(19) corresponding to the $\xi$-th interval is labeled as $T_{on,k_c,\xi}(t)$ and $T_{off,k_c,\xi}(t)$. $\rho_\zeta$ denotes the probability that the time instant $t$ belongs to the $\xi$-th interval and is calculated as follows:

$$\rho_\zeta(t) = (1 - F_\zeta^H(t)) \times F_\zeta^L(t) \quad \forall \zeta = 1,2,3 \cdots \Xi \tag{20}$$

where $F_\zeta^H(t)$ is the cumulative distribution function (CDF) of the higher endpoint of the $\xi$-th interval; $F_\zeta^L(t)$ is the CDF of the lower endpoint of the $\xi$-th interval.

In such a stochastic situation, the aggregate power $P_g(t)$ of TCLs in group $g$ calculated in (6) is rewritten as follows:

$$P_g(t) = \sum_{c=1}^{Q} \left( \left( \sum_{\xi=1}^{\Xi} \frac{T_{on,k_c,\xi}(t)}{T_{on,k_c,\xi}(t) + T_{off,k_c,\xi}(t)} \cdot \rho_\xi(t) \right) \cdot \sum_{\upsilon \in S_c} p_\upsilon \right) \tag{21}$$

The fraction of on-time in an on-and-off cycle is defined as the duty cycle. The duty cycle $\eta_{k_c,\xi}$ of the $c$-th cluster in the $\xi$-th interval is:

$$\eta_{k_c,\xi}(t) = \frac{T_{on,k_c,\xi}(t)}{T_{on,k_c,\xi}(t) + T_{off,k_c,\xi}(t)} \tag{22}$$

In this way, the parameters of Eq. (6) influenced by $\theta_a(t)$ and $\theta_{set,k_c}(t)$ are covered in $\eta_{k_c,\xi}$, while the other parameters are constant values. Considering that the value of $\Delta\theta_a(t)$ and $\Delta\theta_{set,k_c}(t)$ are far lower than the predicted ambient temperature $\bar{\theta}_a(t)$ and the expected set point temperature $\bar{\theta}_{set,k_c}(t)$, the total deviation $\Delta P_g(t)$ of TCLs' aggregate power from the value calculated by Eq. (6) is represented as:

$$\Delta P_g(t) = \sum_{c=1}^{Q} \left( \left( \sum_{\xi=1}^{\Xi} \left( \frac{\partial \eta_{k_c,\xi}(t)}{\partial \theta_a(t)} \bigg|_{\bar{\theta}_a(t),\bar{\theta}_{set,k_c}(t)} \Delta\theta_a(t) + \frac{\partial \eta_{k_c,\xi}(t)}{\partial \theta_{set,k_c}(t)} \bigg|_{\bar{\theta}_a(t),\bar{\theta}_{set,k_c}(t)} \Delta\theta_{set,k_c}(t) \right) \cdot \rho_\xi(t) \right) \cdot \sum_{\upsilon \in S_c} p_\upsilon \right) \tag{23}$$

Since $\Delta\theta_a(t)$ and $\Delta\theta_{set,k_c}(t)$ are independent, the $v$-th order cumulant $\kappa_{\Delta P_g,v}(t)$ of $\Delta P_g(t)$ is given by:

$$\kappa_{\Delta P_g,v}(t) = \left( \sum_{c=1}^{Q} \left( \sum_{\varsigma=1}^{\Xi} \frac{\partial \eta_{k_c,\xi}(t)}{\partial \theta_a(t)} \bigg|_{\bar{\theta}_a(t),\bar{\theta}_{set,k_c}(t)} \cdot \rho_\xi(t) \cdot \sum_{\upsilon \in S_c} p_\upsilon \right) \right)^v \kappa_{\Delta\theta_a,v}$$
$$+ \left( \sum_{c=1}^{Q} \left( \sum_{\varsigma=1}^{\Xi} \frac{\partial \eta_{k_c,\xi}(t)}{\partial \theta_{set,k_c}(t)} \bigg|_{\bar{\theta}_a(t),\bar{\theta}_{set,k_c}(t)} \cdot \rho_\xi(t) \cdot \sum_{\upsilon \in S_c} p_\upsilon \right) \right)^v \kappa_{\Delta\theta_{set},v} \tag{24}$$

Then probability distribution function (PDF) $f_{\Delta P_g(t)}(x)$ of $\Delta P_g(t)$ can be obtained by Gram-Charlier Type A Expansion [46]:

$$f_{\Delta P_g(t)}(x) = \sum_{i=0}^{n} (-1)^i \cdot \frac{c_i}{i!} \cdot H_i(x) \cdot \varphi(x) \tag{25}$$

where $\varphi(x)$ is the CDF of the standard normal distribution, $H_i(x)$ is the Hermite polynomial. $c_i$ is the constant coefficient and can be calculated by the 1-*st* to *n*-th order cumulants of $\Delta P_g(t)$.

The mean value $\bar{P}_g(t)$ of TCLs' aggregate power $P_g(t)$ is obtained by replacing $\theta_a(t)$ and $\theta_{set,k_c}(t)$ a in Eq. (21) as the mean value $\bar{\theta}_a(t)$ and $\bar{\theta}_{set,k_c}(t)$, respectively. Then, CDF $F_{P_g(t)}(x)$ of $P_g(t)$ at each time instant is obtained according to the PDF in (25) and is represented by:

$$F_{P_g(t)}(x) = \int_{-\infty}^{x} f_{\Delta P_g(t)}(x - \bar{P}_g(t))\, dx \tag{26}$$

*4.2. Multi-State Reliability Model of Operating Reserve Provided by TCLs*

The variation of TCLs' aggregate power calculated by Eqs. (21)-(26) is illustrated in Fig. 4(a). $P_i(t)$ denotes the aggregate power of TCLs at bus *i*. The mean value of $P_i(t)$ is represented in black line, while the corresponding variation upper boundary and lower boundary of $P_i(t)$ is represented by red and blue line, respectively. Fig. 4(a) illustrates that the provision of the equivalent operating reserve is a continuous process, which involves the gradual reduction of aggregate power and the gradual rebound of aggregate power. Therefore, instead of the conventional two-state model for operating reserve, the equivalent operating reserve provided by TCLs of bus *i* at time *t* are represented as a multi-state operating reserve provider, which is abbreviated as MORT$_i$(*t*).

Let $RC_{i,j_i^{TCL}}$, $j_i^{TCL} = 1,...,K_i^{TCL}$ be the state of the reserve capacity at bus *i* for the time *t*. The range of state space is determined according to the maximum reserve capacity during the reserve deployment process. The deviation $\tau_{i,j_i^{TCL}}$ between the neighboring states $RC_{i,j_i^{TCL}}$ and $RC_{i,j_i^{TCL}+1}$ equals to the standard deviation corresponding to the state $RC_{i,j_i^{TCL}}$, which can be obtained by the PDF shown Fig. 4(b). Such division of state space according to standard deviation has been proved to give a satisfactory computational accuracy compared to the reliability computational result obtained by Monte-Carlo method [47].

Let $\rho_{i,j_i^{TCL}}(t)$, $j_i^{TCL} = 1,...,K_i^{TCL}$ be the probabilities of the reserve capacity $RC_{i,j_i^{TCL}}$ at bus *i* for the time *t*. The probability of each state for the time *t* can be obtained from the CDF of $P_i(t)$ obtained in Eq. (26). As illustrated in Fig. 4(b), the CDF of TCLs' aggregate power within the variation boundary between $P_i^{max}(t_1)$ and $P_i^{min}(t_1)$ at the time $t_1$ is $F_{P_i(t_1)}(x)$. Then, the probability of reserve capacity at the state $j_i^{TCL}$ is calculated by the difference in CDF corresponding to the endpoints of the state, as illustrated in Fig. 4(c). Therefore, $\rho_{i,j_i^{TCL}}(t)$ can be obtained by the CDF of aggregate power generated from Eq. (26) and represented by:

$$\rho_{i,j_i^{TCL}}(t) = F_{P_i(t)}(P_i^0 - RC_{i,j_i^{TCL}} + \tau_{i,j_i^{TCL}}) - F_{P_i(t)}(P_i^0 - RC_{i,j_i^{TCL}}) \tag{27}$$

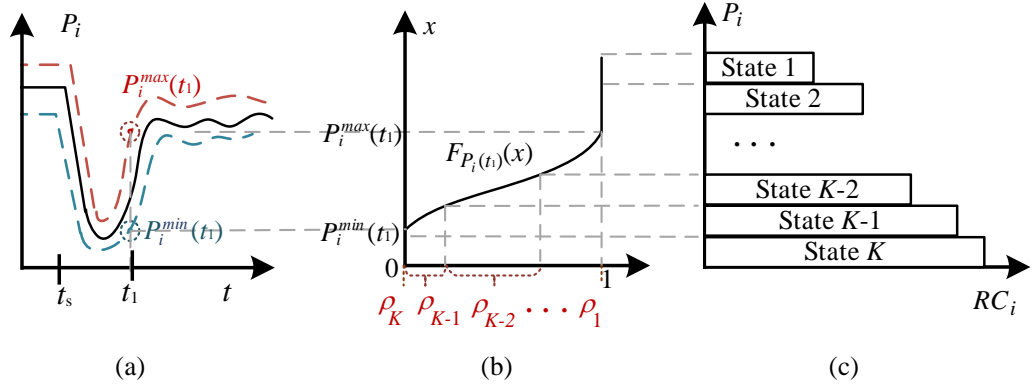

Fig. 4. Probabilistic Analysis of ORT Considering TCLs' Dynamic Response (a) Aggregate power of TCLs after changing set point temperature at $t_s$ (b) CDF of TCLs' aggregate power at the time $t_1$ (c) Obtaining the probabilities corresponding to the states of reserve capacity at the time $t_1$

$L_Z$-transform, which is an extension of traditional universal generating function (UGF) technique [17], [18], has been proved to effectively represent multi-state units for discrete-state continuous-time reliability evaluation [20], [21]. Hence, $L_Z$-transform is applied in this paper to represent the power output distribution of $\text{MORT}_i(t)$ and can be defined as the following polynomial:

$$Lz_i^{MORT}(z,t) = \sum_{j_i^{TCL}=1}^{K_i^{TCL}} \rho_{i,j_i^{TCL}}(t) \cdot z^{RC_{i,j_i^{TCL}}} \tag{28}$$

where $z$ in general is a complex variable. The introduction of this complex variable provides a comprehensive approach for the system state enumeration that can substitute complicated combinational algorithms [20].

## 5. Power System Short-term Reliability Evaluation Considering Hybrid Generation and Operating Reserve Providers

### 5.1. System Reliability Model Considering Hybrid Generation and Operating Reserve Providers

The total available generation capacity of the power systems comes from hybrid generation providers (e.g., conventional generating units, wind farms, etc.) and hybrid operating reserve providers (e.g., conventional operating reserve providers, ORT, etc.) [48]. In this paper, the multi-state model of hybrid generation providers and operating reserve providers are combined as an equivalent power generation provider, which is obtained using $L_Z$-transform.

There exists the co-operation between ORT and conventional operating reserve providers. It takes a period of time, defined as lead time, for reserve service providers to start providing operating reserve after the reserve deployment instruction [47]. Usually, the lead time of operating reserve provided by TCLs is much shorter than that of conventional generation units. The operating reserve provided by TCLs is deployed the earliest to enhance system reliability in a short time. Then, conventional generation units are deployed to supplement the operating reserve provided by TCLs. In this way, the number of reserve service provider $N_i(t)$ at bus $i$ for the time $t$ is determined by the deployment procedure. Similar with the $L_Z$-transform to represent ORT, the distribution of reserve capacity corresponding to the $n$-th conventional operating reserve provider is obtained in the previous work [10] and represented as follows:

$$Lz_{i,n}^{MRS}(z,t) = \sum_{j_{i,n}^{RS}=1}^{K_{i,n}^{RS}} \rho_{i,j_{i,n}^{RS}}(t) \cdot z^{RC_{i,j_{i,n}^{RS}}} \tag{29}$$

where $K_{i,n}^{RS}$ is the total number of states corresponding to the $n$-th reserve service provider. $\rho_{i,j_{i,n}^{RS}}(t)$, $j_{i,n}^{RS} = 1,...,K_{i,n}^{RS}$ are the state probabilities of the reserve capacity $RC_{i,j_{i,n}^{RS}}$ corresponding to the $n$-th operating reserve provider at time $t$.

The total reserve capacity of hybrid operating reserve providers is the accumulation of each provider's reserve capacity. Therefore, the $L_Z$-transform of hybrid operating reserve providers $\text{MHOR}_i(t)$ are calculated using the parallel composition operator $\Omega_{\phi p}$ and represented as follows:

$$\begin{aligned}
Lz_i^{MHOR}(z,t) &= \Omega_{\phi p}\left\{Lz_i^{MORT}(z,t), Lz_{i,1}^{MRS}(z,t) \cdots Lz_{i,n}^{MRS}(z,t)\right\} \\
&= \Omega_{\phi p}\left\{\sum_{j_i^{TCL}=1}^{K_i^{TCL}} \rho_{i,j_i^{TCL}}(t) \cdot z^{RC_{i,j_i^{TCL}}}, \sum_{j_{i,1}^{RS}=1}^{K_{i,1}^{RS}} \rho_{i,j_{i,1}^{RS}}(t) \cdot z^{RC_{i,j_{i,1}^{RS}}}, ..., \sum_{j_{i,N(t)}^{RS}=1}^{K_{i,N(t)}^{RS}} \rho_{i,j_{i,N(t)}^{RS}}(t) \cdot z^{RC_{i,j_{i,N(t)}^{RS}}}\right\} \\
&= \sum_{j_i^{TCL}=1}^{K_i^{TCL}} \sum_{j_{i,1}^{RS}=1}^{K_{i,1}^{RS}} ... \sum_{j_{i,N(t)}^{RS}=1}^{K_{i,N(t)}^{RS}} \prod_{n=1}^{N(t)} \rho_{i,j_i^{TCL}}(t) \cdot \rho_{i,j_{i,n}^{RS}}(t) \cdot z^{\phi_P(RC_{i,j_i^{TCL}}, RC_{i,j_{i,1}^{RS}}, ..., RC_{i,j_{i,N(t)}^{RS}})} \\
&= \sum_{j_i^{MHOR}=1}^{K_i^{MHOR}} \rho_{i,j_i^{MHOR}}(t) \cdot z^{HRC_{i,j_i^{MHOR}}}
\end{aligned} \tag{30}$$

where $K_i^{MHOR}$ is the total number of states corresponding to the hybrid operating reserve at bus $i$. $\rho_{i,j_i^{MHOR}}(t)$, $j_i^{MHOR} = 1,...,K_i^{MHOR}$ are the state probabilities of the reserve capacity $HRC_{i,j_i^{MHOR}}$ corresponding to the hybrid operating reserve at bus $i$ for the time $t$.

The $L_Z$-transform $Lz_i^{MHGU}(z,t)$ for the multi-state hybrid generation provider $\text{MHGU}_i(t)$ at bus $i$ for time $t$, including the wind farms and the conventional generating units, have been put forward in the previous work [10]. The combination of the $\text{MHGU}_i(t)$ and the $\text{MHOR}_i(t)$ can be represented as a multi-state hybrid generation and reserve provider $\text{MHGR}_i(t)$ at bus $i$ for time $t$. The $L_Z$-transform for the $\text{MHGR}_i(t)$ is obtained using the parallel composition operator $\Omega_{\phi p}$ and represented as follows [10]:

$$Lz_i^{MHGR}(z,t) = \Omega_{\phi p}\left\{Lz_i^{MHGU}(z,t), Lz_i^{MHOR}(z,t)\right\}$$

$$= \Omega_{\phi p}\left\{\sum_{j_i^G=1}^{K_i^G} \rho_{i,j_i^G}(t) \cdot z^{AG_{i,j_i^G}}, \sum_{j_i^{MHOR}=1}^{K_i^{MHOR}} \rho_{i,j_i^{MHOR}}(t) \cdot z^{HRC_{i,j_i^{MHOR}}}\right\}$$

$$= \sum_{j_i^G=1}^{K_i^G} \sum_{j_i^{MHOR}=1}^{K_i^{MHOR}} \rho_{i,j_i^G} \cdot \rho_{i,j_i^{MHOR}}(t) \cdot z^{\phi_P(AG_{i,j_i^G}, HRC_{i,j_i^{MHOR}})} \quad (31)$$

$$= \sum_{j_i=1}^{K_i} \rho_{i,j_i}(t) \cdot z^{AG^*_{i,j_i}}$$

where $K_i^G$ is the total number of states corresponding to the multi-state hybrid generation provider MHGU$_i$(t) at bus $i$. $\rho_{i,j_i^G}(t)$, $j_i^G = 1,...,K_i^G$ are the state probabilities corresponding to the available generation capacity $AG_{i,j_i^G}$ of the multi-state hybrid generation provider. The accumulation of available generation capacity from generation providers and operating reserve providers can be regarded as equivalent available generation capacity $AG^*_{i,j_i}$. $K_i$ is the total number of states corresponding to the equivalent available generation capacity at bus $i$. $\rho_{i,j_i}(t)$, $j_i = 1,...,K_i$ are the state probabilities of $AG^*_{i,j_i}$ corresponding to the hybrid generation and operating reserve provider at bus $i$ for the time $t$.

*5.2. Reliability Indices*

After obtaining the $L_Z$-transform for the MHGR$_i$(t), the load curtailment at each bus is calculated by optimal power flow composition operator $\Omega_{\Phi OPF}$ [10]. For an $N$-bus system with $K$ system states, the $L_Z$-transform to obtain the load curtailment at bus $i$ for the time $t$ is represented as follows:

$$Lz_i^{MLC}(z,t) = \Omega_{\phi OPF}\left\{Lz_{i,1}^{MHGR}(z,t),...,Lz_{i,N}^{MHGR}(z,t)\right\}$$

$$= \sum_{i=1}^{N}\sum_{j_i=1}^{K_i}\sum_{j_L=1}^{K_L} \rho_{i,j_i}(t) \cdot \rho_{j_L}(t) \cdot z^{\phi_{OPF}\{AG^*_{1,j_1},...,AG^*_{N,j_N}\}} \quad (32)$$

$$= \sum_{j=1}^{K} \rho_j(t) \cdot z^{LC_{j_i}(t)}$$

where $\rho_j(t)$ and $LC_{j_i}(t)$ are the probability and load curtailment at bus $i$ for the system state $j$ at time $t$, respectively. $K_L$ is the number of states for the transmission network; $\rho_{j_L}(t)$ is the probability of the transmission network state $j_L$ at time $t$.

The optimal power flow composition operator $\Omega_{\Phi OPF}$ used in Eq. (32) is utilized to minimize the total system load curtailment for the system state $j$ at time $t$:

$$\text{Min } f_j = \sum_{i=1}^{N} LC_{j_i}(t) \quad (33)$$

$$\text{s.t.} \quad \mathbf{B_j \theta_j}(t) = \mathbf{P_j}(t) - \mathbf{D_j}(t) \quad (34)$$

$$\mathbf{LC_j}(t) = \mathbf{\bar{D}}(t) - \mathbf{D_j}(t) \quad (35)$$

$$0 \leq p_{j_i}(t) \leq AG^*_{i,j_i} \quad (36)$$

$$\left| \frac{1}{x_{j_{ik}}} \left( \theta_{j_i}(t) - \theta_{j_k}(t) \right) \right| \leq \left| F_{ik}^{\max} \right| \tag{37}$$

where Eq. (34) is the DC power flow constraints, Eq. (35) is the load curtailment constraints, Eq. (36) is the generation output limits:, Eq. (37) is the line flow constraints; $\mathbf{B_j}$ is the admittance matrix of transmission network, $\mathbf{\theta_j}(t)$ is phase angle vector of bus voltages at time $t$, $\mathbf{P_j}(t) = [p_{j_1}(t),...,p_{j_N}(t)]^T$ is the vector of equivalent power generation for the state $j$ at the time $t$, $\mathbf{D_j}(t) = [D_{1,j}(t),...,D_{N,j}(t)]^T$ and $\mathbf{\bar{D}}(t) = [\bar{D}_1(t),...,\bar{D}_N(t)]^T$ represent the vector of the bus loads for the state $j$ at time $t$ and the vector of the bus loads for the normal state for the time $t$, respectively. $\mathbf{LC_j}(t) = [LC_{1,j}(t),...,LC_{N,j}(t)]^T$ is the vector of load curtailment for the state $j$ for the time $t$. $p_{j_i}(t)$ is power generation of the $MHGR_i(t)$ and $\theta_{j_i}(t)$ is the phase angle of voltage at bus $i$ for the time $t$, $x_{j_{ik}}$ and $\left| F_{ik}^{\max} \right|$ are the reactance and maximum power flow of the line between buses $i$ and $k$ respectively.

The system reliability indices defined in [10], including the LOLP, EENS and LOLE, are utilized to evaluate system reliability.

$LOLP_i(t)$ is defined as the loss of load probability at bus $i$ for time $t$, which can be evaluated as:

$$LOLP_i(t) = \sum_{j=1}^{K} \rho_j(t) \mathbf{1}(LC_{j_i}(t) > 0) \tag{38}$$

where $\mathbf{1}(True) \equiv 1, \ \mathbf{1}(False) \equiv 0$.

$EENS_i(\tau)$ is defined as the expected energy not supplied at bus $i$ during the operation period $\tau$, which can be evaluated as:

$$EENS_i(\tau) = \int_0^{\tau} \left( \sum_{j=1}^{K} \rho_j(t) \cdot LC_{j_i}(t) \right) \cdot dt \tag{39}$$

$LOLE_i(\tau)$ is defined as the loss of load expectation at bus $i$ during the operation period $\tau$, which can be evaluated as:

$$LOLE_i(\tau) = \int_0^{\tau} \left( \sum_{j=1}^{K} \rho_j(t) \cdot \mathbf{1}(LC_{j_i}(t) > 0) \right) \cdot dt \tag{40}$$

*5.3. Computation Procedure for Reliability Evaluation*

The basic procedures for the time varying reliability assessment of power systems are as follows:

**Step1**: Input characteristic parameters, including the reserve deployment time instant, number of TCLs, distribution of set point temperature, distribution of ambient temperature, etc.

**Step2:** Obtain the probability distribution of TCLs' consumed power after the changes of set point temperature using Eqs. (6)-(26).

**Step3:** Determine the $L_Z$-transform for the operating reserve $MORT_i(t)$ provided by TCLs considering the stochastic consumers' behavior and ambient temperature using Eqs. (27)-(28).

**Step4:** Determine the system $L_Z$-transform for the hybrid generation and operating reserve providers using Eqs. (29)-(31).

**Step5:** Obtain the $L_Z$-transform for determining the load curtailment at each bus using Eqs.(32)-(37).

**Step6:** Calculate the $LOLP_i(t)$, $EENS_i(\tau)$, $LOLE_i(\tau)$ using Eqs. (38)-(40), respectively.

## 6. Case studies

Four case studies are conducted to illustrate the proposed multi-state operating reserve model of aggregate TCLs for power system short-term reliability evaluation. Firstly, the aggregate power of TCLs and the corresponding equivalent operating reserve are obtained by the proposed analytical model, so that the particular dynamic characteristics of operating reserve provided by TCLs can be observed. Secondly, considering the variations of ambient temperature and set point temperature, states of reserve capacity and the corresponding time-varying probabilities are obtained by the proposed multi-state operating reserve model of TCLs. Thirdly, short-term reliability of IEEE Reliability Test System (RTS) is evaluated considering operating reserve provided by TCLs. In this case, the impact of TCLs' dynamic characteristics on the system reliability can be obtained. Finally, short-term reliability of power system in Nantong with a typical summer day is evaluated, so that the real application of the proposed method in enhancing the system reliability treated by peak demand and assisting operating reserve commitment decisions are illustrated.

### 6.1. Dynamics of Operating Reserve Provided by TCLs

This case illustrates the dynamics of TCLs' aggregate power and the corresponding equivalent operating reserve. Parameters of TCLs are set according to [49] and presented in Table 1. The total number of controllable TCLs for the provision of operating reserve is set as 100,000. The number of cluster $Q$ in Eq. (6) is set as 8 according to the Calinski-Harabasz criterion [42]. All the controllable TCLs are controlled for providing operating reserve through increasing the set point temperature by 1 ºC. In this way, the aggregate power of TCLs will decrease and therefore provide equivalent operating reserve. During this process, TCLs' aggregate power and the corresponding equivalent operating reserve are plot by the curves and stack area in Fig. 5, respectively.

Fig. 5 illustrates that the initial aggregate power of TCLs is approximately 180MW. After the provision of operating reserve at 1:00, aggregate power of TCLs gradually decreases and reaches the minimum point at around 1:20. During this process, the equivalent operating reserve gradually increases to the level of 180MW. However, the aggregate power of TCLs rebounds after 1:20. Such phenomenon is referred to as demand response rebound in many literatures and cannot be reflected by traditional

Table. 1 TCL Physical Parameters

| Parameters | Descriptions | Values | Units |
| --- | --- | --- | --- |
| $C_\upsilon$ | Thermal capacity | $U$ (1.5, 2.5) | kWh/ºC |
| $R_\upsilon$ | Thermal resistance | $U$ (1.5, 2.5) | ºC/kW |
| $COP_\upsilon$ | Coefficient of performance | 2.5 | / |
| $p_\upsilon$ | Input power | $U$ (4, 7.2) | W |
| $\theta_a$ | Ambient temperature | $N$ (32, 1) | ºC |
| $\theta_{set,\upsilon}$ | Set point temperature | $U$ (18, 27) | ºC |

Normal distribution with the mean value of $\mu$ and the standard deviation of $\sigma$ is abbreviated to $N(\mu,\sigma)$; uniform distribution with the minimum and maximum value of $a$ and $b$, respectively, is abbreviated to $U(a,b)$.

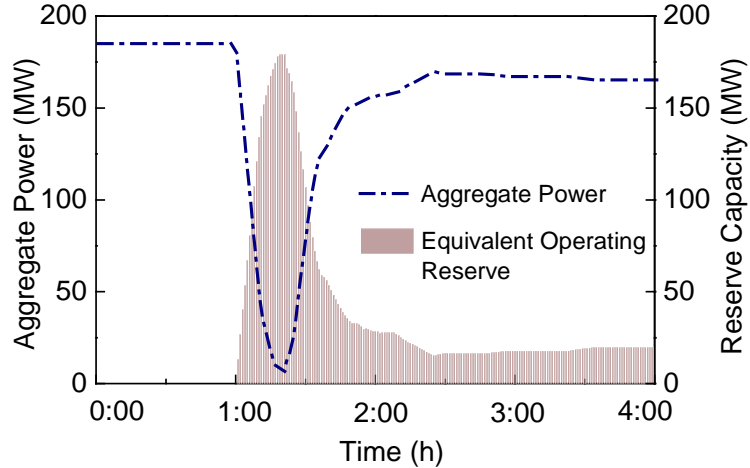

Fig. 5 TCLs' aggregate power and the corresponding equivalent operating reserve

methods based on derated rates [6]. Because of the rebound effect, the equivalent operating reserve also decreases after 1:20 and eventually reaches the value of only 20MW. Hence, the characteristics of operating reserve provided by TCLs are different from that provided by conventional generating units. It is essential to involve the dynamic response of TCLs in the short-term reliability so that the system reliability level can be accurately evaluated.

*6.2. Illustration of Multi-state Operating Reserve Provided by TCLs*

This case illustrates the multi-state operating reserve model of TCLs considering variations of ambient temperature and set point temperature. Corresponding to the mean value of operating reserve reflected by the stack area in Fig. 5, the multi-state operating reserve model of TCLs obtained by Eqs. (20)-(28) is reflected in Fig. 6. The y-axis is the states of reserve capacity. The colormap interprets the probabilities corresponding to the states at each time instant.

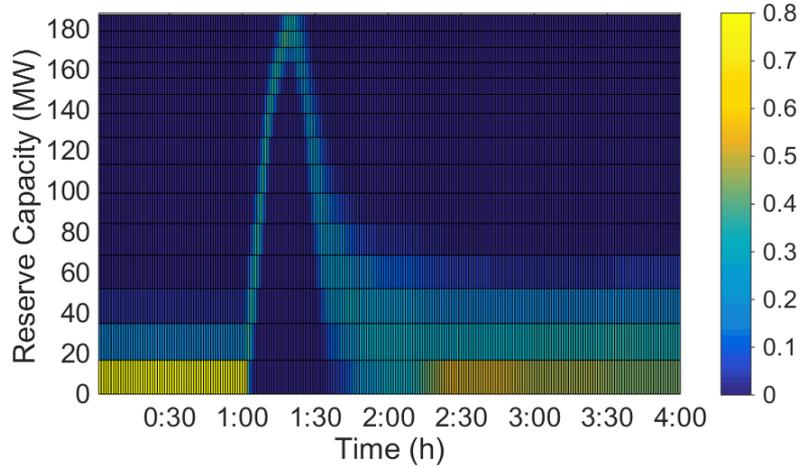

Fig. 6 Multi-state operating reserve provided by TCLs

It can be observed from the y-axis of Fig. 6 that there are fifteen states of reserve capacity distribute between 0MW to 180MW. Before 1:00, the probability of reserve capacity at 0MW is around 0.8, which is corresponding to the situation without reserve deployment. After the control of TCLs for providing operating reserve at 1:00, the states with higher probability shown by the area with lighter color gradually increase. At approximately 1:20, the states with the highest reserve capacity ranging from 160MW and 180MW share the highest probability. This means that TCLs at this time instant are expected to provide the highest reserve capacity. After that, the area with lighter color gradually moves to the states with lower reserve capacity, which corresponds to the decrease of reserve capacity resulted from demand response rebound. Therefore, the proposed multi-state operating reserve model can represent the changes in probability of each reserve capacity over time resulting from the aggregate dynamics of TCLs.

*6.3. Power System Short-Term Reliability Evaluation Considering Operating Reserve Provided by TCLs*

This case evaluates the power system short-term reliability with operating reserve provided by TCLs. The modified IEEE Reliability Test System in [10] is utilized to illustrate the proposed models and techniques. The total demand is 2850 MW. The total number of controllable TCLs for the provision of operating reserve is set as 100,000. In this way, the aggregate power of all the controllable TCLs is around 200MW, which equals to 7% of total demand. The distribution of controllable TCLs at each bus is proportional to the base load in these buses. Apart from TCLs, there are five 40-MW gas thermal generators working as operating reserve provider [10], which are located at bus 1 (three units) and bus 2 (two units), respectively. Hybrid generation providers consist of conventional generators and wind farms. Reliability model of these generation units are the same as that proposed in [10]. A 500-MW wind farm including 250 identical 2-MW wind turbines is added to Bus 21. The cut-in, rated, and cut-out wind speeds of a wind turbine are 4, 15 and 25 km/h, respectively. The online conventional generators consist of four 576-MW coal thermal generators and three 197-MW oil thermal generators. The four coal thermal generators are located at buses 15, 16, 18 and 23. The three oil thermal generators are installed at bus 13.

Two senarios are included: 1) WoOR: the base scenario without the commitment of operating reserves. 2) ORT: TCLs are controlled for providing operating reserve at the time 1:00. In the senond senario, reliability of TCLs obtained by Monte-Carlo method (MC) and the proposed analytical method are compared, labeling as ORT-MC, ORT-ANL, repectively. The initial wind speed at the wind farm is set as 16 km/h. Correspondingly, wind turbines generate rated power at time t=0. All the other generating units are in good condition at the beginning of the operating time. The LOLP for a representative load bus (i.e., bus 6) from 0:00 to 4:00 is illustrated in Fig. 7. EENS and LOLE corresponding to the scenarios are illustrated in Table. 2, where the value of EENS and LOLE are the value at the end of simulation interval.

The simulation scenarios are conducted on a PC with Intel 2.3 GHz 2-core processor (4MB L3 cache), 8 GB memory. The computational time of reliability evaluation obtained by MC and the proposed analytical method are 3271 seconds and 1918 seconds, respectively. Therefore, the computational time of MC is much longer than the proposed method.

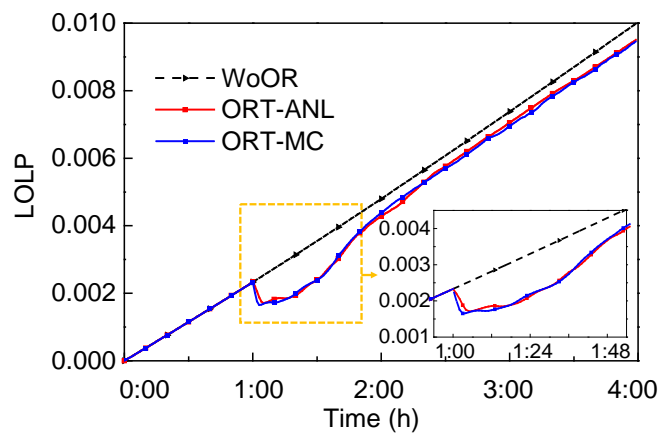

Fig. 7 Instant LOLP at bus 6 from 0 to 4h in IEEE Reliability Test System

Table. 2 Reliability indices obtained by different methods in IEEE Reliability Test System

| Method | EENS (MWh) | LOLE (h) |
| --- | --- | --- |
| ORT-ANL | 0.03538 | 0.007362 |
| ORT-MC | 0.03520 | 0.007274 |
| WoOR | 0.04214 | 0.009151 |

EENS and LOLE are the value at the end of simulation interval.

The reliabillity indices obtained from MC is regarded as the benchmark to evaluate the accuracy of reliabillity indices obtained by the proposed method. The curve of ORT-ANL and ORT-MC in Fig. 7 illustrates that the LOLP calculated by the proposed method is close to that calculated by MC. Table. 2 illustrates that EENS and LOLE calculated by the proposed analytical method is 0.03538 MWh and 0.007362h, which is highly close to the EENS (0.03520MWh) and LOLE (0.007274h) calculated by MC. Therefore, the proposed method can reduce the computational time and characterize the dynamics of TCLs at the same time, which guarantee that the impact of operating reserve provided by TCLs on the short-term reliability can be adequately evaluated.

The curve of ORT-MC in Fig. 7 illustrates that the instant LOLP decreases from 0.0024 to 0.0015 after the commitment of operating reserve provided by TCLs. The LOLP after the commitment of ORT remains approximately 0.001 lower than the initial LOLP (shown by the curve of WoOR in Fig. 7) without reserve commitment during the period between 1:00 to 1:25. However, because of the demand response rebound, the instant LOLP bagins to increase at 1:30. Meanwhile, Table. 2 illustrates that EENS and LOLE with ORT is 0.03520 MWh and 0.007274 h, which is 0.00694 MWh (=0.04214 MWh -0.03520 MWh) and 0.001921h (=0.009151h -0.007274h) lower than the EENS (0.04214 MWh) and LOLE (0.009151h) without reserve commitment. This demonstrates that the commitment of ORT can enhance the system reliability, despite of the negative effect brought by the demand response rebound.

*6.4. Usage of the Proposed Method for the Real Applications*

This subsection illustrates the usage of the proposed method for the real applications based on the power system of Nantong, a large city in Jiangsu Province. Firstly, the proposed method is applied to the power system of Nantong on a summer day with high electricity demand, where the operating reserve provided by TCLs is committed to enhance the system reliability treated by peak electricity demand. Secondly, the proposed method is applied to assist system operators in adequately evaluating system operating pressures during a short interval and cooperating different types of operating reserve.

*6.4.1 Enhancing power system reliability threatened by peak electricity demand using operating reserve provided by TCLs*

The total capacity of conventional generators in power systems of Nantong is 10.851GW and the capacity of wind power generation is 2.079 GW [50]. This power system consists of 32 buses and 96 transmission lines [51]. The experiment is conducted on a representative summer day when the peak demand is 7.5GW with the load curve illustrated in Fig. 8. The total amount of controllable load is 0.87GW [52]. Half the controllable loads (0.44GW) are assumed to consist of TCLs, which accounts for approximate 5.9% of peak demand and are provided by approximately 240,000 TCLs. The electricity demand in Fig. 8 illustrates that electricity consumption increases sharply since 6:00 and reaches approximately the peak value at 10:00. Therefore, this case simulates the power system short-term reliability between 8:00-12:00 to interpret the impact of high electricity demand and therefore validate the effectiveness of operating reserve provided by TCLs to relieve the power system reliability treats.

TCLs are controlled for providing operating reserve at the time 9:00. The lables of senarios are the same as that in subsection 6.3. The initial wind speed at the wind farm is set as 13 km/h. Correspondingly, wind turbines generate power at derated states lower than the rated level. All the other generating units are in good condition at the beginning of the operating time. The LOLP for a representative load bus from 8:00 to 12:00 is illustrated in Fig. 9. EENS and LOLE corresponding to the scenarios are illustrated in Table. 3, where the value of EENS and LOLE are the value at the end of simulation interval.

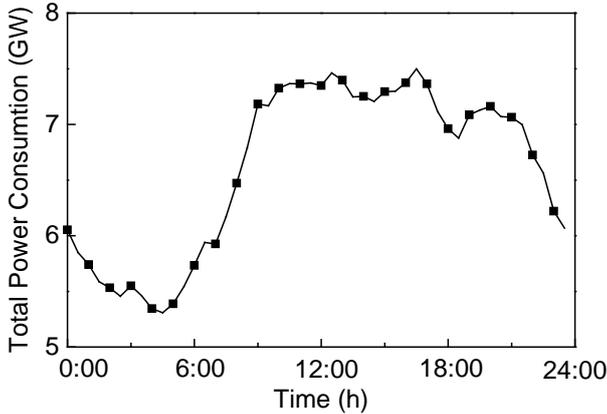

Fig. 8 Total power consumption in Nantong on a typical summer day

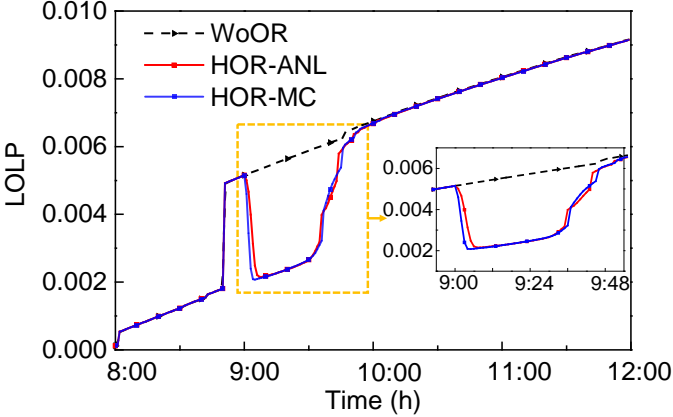

Fig. 9 Instant LOLP for a representative bus from 8:00 to 12:00 with operating reserve provided by TCLs deployed at 9:00

Table. 3 Reliability indices obtained by different methods

| Method | EENS (MWh) | LOLE (h) |
|---|---|---|
| HOR-ANL | 0.8209 | 0.002704 |
| HOR-MC | 0.8171 | 0.002683 |
| WoOR | 0.9982 | 0.003218 |

EENS and LOLE are the value at the end of simulation interval

As illustrated by the curve of WoOR in Fig. 9, the increase of demand leads to the giant increase of LOLP from 0.002 to 0.005 at approximately 8:55. Compared with the curve of WoOR, the curve of ORT-ANL and ORT-MC shows that the deployment of operaitng reserve provided by TCLs at 10:00 successfully reduces LOLP to 0.002, which is the level without the peak demand. During the period between 9:00 to 9:30, LOLP remain around 0.003 lower than the curve without reserve deployemnt. However, because of the demand response rebound, the instant LOLP bagins to increase at 9:30, and eventually reaches the level without the unit commitment. Table. 3 illustrates that EENS and LOLE calculated by the proposed analytical method is 0.8209MWh and 0.002704h, which is highly close to the EENS (0.8171MWh) and LOLE (0.002683 h) calculated by MC. In comparison, EENS and LOLE without reserve deployment is 0.9982 MWh and 0.003218 h, which means that the reserve deployment in HOR-ANL and HOR-MC reduce the reliability indices by approximate 18%. Therefore, the deployment of operating reserve provided by TCLs can enhance the system reliability threatened by peak demand.

6.4 *Assisting system operators in operating reserve commitment decision*

This case illustrates the function of the proposed method to assist system operators in operating reserve commitment decision. To guarantee the reliable operation of power systems, system operators can cooperate the commitment of different kinds of operating reserves and therefore mitigate the demand response rebound of TCLs. Operating reserve provided by TCLs is still deployed at 9:00. The other operating reserve is provided by gas thermal generators and is committed at 9:30. The *LOLP* for a representative load bus from 8:00 to 12:00 considering the commitment of hybrid operating reserve is illustrated in Fig. 10. Reliability of TCLs obtained by Monte-Carlo method (MC) and the proposed analytical method are labeled as HOR-MC, HOR-ANL, repectively. EENS and LOLE corresponding to this case are illustrated in Table. 4, where the value of EENS and LOLE are the value at the end of simulation interval.

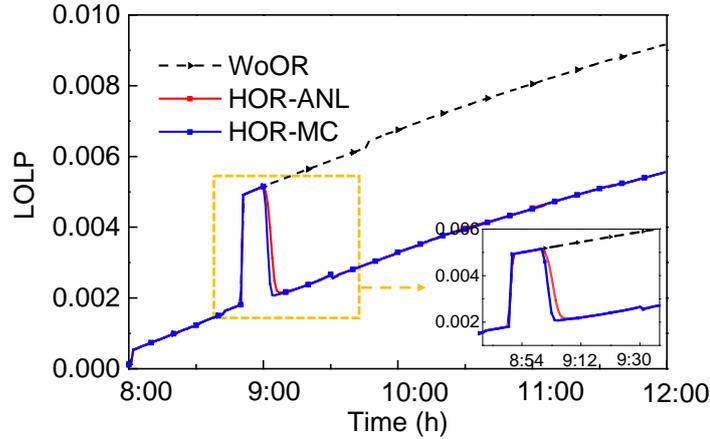

Fig. 10 Instant LOLP at a representative bus from 8:00 to 12:00 considering hybrid operating reserve

Table. 4 Reliability indices obtained by different methods

| Method | EENS (MWh) | LOLE (h) |
| --- | --- | --- |
| HOR-ANL | 0.2469 | 0.0006910 |
| HOR-MC | 0.2449 | 0.0006789 |
| WoOR | 0.9982 | 0.003218 |

EENS and LOLE are the value at the end of simulation interval

Similar to Fig. 9, LOLP after the reserve commitment decreases from 0.005 to 0.002 at 9:00 in this case, illustrated by the curve of HOR-ANL and HOR-MC in Fig. 10. However, compared with the large increase of LOLP derived from the demand response rebound in Fig. 9, the LOLP in Fig. 10 remains at the reduced level after 9:30. Moreover, Table. 4 illustrates that EENS and LOLE corresponding to Fig. 10 are reduced to the level of approximately 0.24MWh and 0.00068 h, which is much lower than the value shown in Table. 3 (approximately 0.82 MWh and 0.0027 h). This is because the commitment of operating reserve provided by conventional generation units compensates for the rebound capacity of TCLs. Compared with the EENS (0.9982 MWh) and LOLE (0.003218h) without operating reserve commitment, EENS and LOLE in this scenario are reduced by approximate 75%. Hence, involving the dynamics of TCLs can reflect the particular effects of operating reserve provided by TCLs to the system reliability, e.g., the demand response rebound. In this way, special deployment strategy, such as the co-operation within hybrid operating reserve providers, can be designed to reduce the negative impact of the demand response rebound.

## 7. Conclusions

This paper presents a novel multi-state reliability model of operating reserve provided by TCLs for the power system short-term reliability evaluation. The dynamic response of TCLs is characterized according to the migration of TCLs' room temperature during the reserve deployment process. On this basis, the probability distribution of operating reserve provided by TCLs is obtained by cumulants. $L_Z$-transform approach is further applied to represent the system reliability with hybrid generation units and operating reserve providers. The accuracy of the proposed method is validated against the Monte Carlo method. Illustrative results demonstrate that the proposed method can effectively model the dynamic characteristics of operating reserve provided by TCLs. Results obtained from the paper can be summarized in the following aspects:

1) Operating reserve provided by TCLs involves a gradual decrease of reserve capacity corresponds to the demand response rebound, which is different from conventional operating reserve.
2) With multi-state operating reserve model of TCLs, the impact of the dynamic characteristics of operating reserve provided by TCLs can be reflected in the power system short-term reliability evaluation with high accuracy.
3) System operators can adequately aware of system operating pressures during a short interval with the proposed power system short-term reliability evaluation technique.
4) Different operating reserves can be dispatched in cooperation to eliminate the impact of demand response rebound from TCLs.

In conclusion, this paper provides a comprehensive approach to analytically obtain the aggregate dynamic response of TCLs for providing operating reserve under uncertainties, and therefore contribute to the accurate evaluation of power system short-term reliability.

## 8. Acknowledgements

The research is supported by the National Natural Science Foundation of China under Grant 51577167 and Grant 51537010.